\documentclass{article}
\usepackage{amsfonts}
\usepackage[letterpaper,body={14cm,22cm}, mag=1000]{geometry}
\usepackage{amssymb}
\usepackage{secdot}
\usepackage{setspace}
\usepackage{titlesec}
\titleformat{\section}[runin]{\bfseries\filcenter}{\thesection}{1em}{}

\renewcommand{\thesection}{\arabic{section}}
\title{\large \bf Finite $p$-groups with central automorphism group of minimal order}
\author{\small \bf Deepak Gumber and Hemant Kalra \\
\small \em School of Mathematics and Computer Applications\\
\small \em Thapar University, Patiala - 147 004,
India}

\date{}
\newtheorem{thm}{Theorem}[section]
\newtheorem{lm}[thm]{Lemma}

\newtheorem{cor}[thm]{Corollary}

\begin{document}

\maketitle
\begin{abstract}
\noindent {\bf Abstract.} We study finite $p$-groups $G$ of coclass upto 4 for which the group $\mathrm{Aut}_z(G)$ of all central automorphisms of $G$ is of minimal possible order. As a consequence, we obtain very short and elementary proofs of main results of Sharma and Gumber \cite{Sharma}.

\end{abstract}
\vspace{2ex}


\vspace{2ex}


\section{Introduction.} Let $G$ be a finite $p$-group.
An automorphism $\alpha$ of $G$ is called a central automorphism if for every $x\in G$, the element $x^{-1}\alpha(x)$ lies in the center $Z(G)$ of $G$. The center $Z(\mathrm{Inn}(G))$ of the group of all inner automorphisms of $G$ is always contained in $\mathrm{Aut}_z(G)$. Curran \cite{Curran} first considered the case when  $\mathrm{Aut}_z(G)$ is of minimal order. He proved that for $\mathrm{Aut}_z(G)$ to be equal to $Z(\mathrm{Inn}(G))$, $Z(G)$ must be contained in the derived group $G'$ and $Z(\mathrm{Inn}(G))$ must not be cyclic. Observe that if $G$ is of nilpotency class 2, then 
$Z(\mathrm{Inn}(G))=\mathrm{Inn}(G)$ and therefore $\mathrm{Aut}_z(G)=Z(\mathrm{Inn}(G))$ if and only if $G'=Z(G)$ and $Z(G)$ is cyclic by \cite[Main theorem]{CurMcc}. Also, if $G$ is of maximal class, then $Z(\mathrm{Inn}(G))$ is cyclic and hence $\mathrm{Aut}_z(G)>Z(\mathrm{Inn}(G))$. Therefore, to characterize finite $p$-groups $G$ for which $\mathrm{Aut}_z(G)$ is of minimal order, we can assume that neither $G$ is of maximal class and nor of class 2. And if this is the case, then of course $|G|\ge p^5$.

In Theorem 2.1, we give necessary and sufficient conditions on a finite $p$-group $G$ with cyclic center such that $\mathrm{Aut}_z(G)$ is of minimal order.  As a consequence of it, we obtain necessary and sufficient conditions on finite $p$-groups $G$ of coclass upto 4 for which $\mathrm{Aut}_z(G)=Z(\mathrm{Inn}(G))$. 
Recently, Sharma and Gumber \cite{Sharma} have characterized such finite $p$-groups of order $p^5$ and $p^6$($p$ an odd prime). We, in particular, also characterized such finite $p$-groups of order upto $p^7$ for any prime $p$. Our proofs are short and elementary than given by Sharma and Gumber \cite{Sharma}.

If $G$ is of order $p^n$ and of nilpotency class $c$, then $G$ is said to be of coclass $n-c$. By $\Phi(G)$ and $Z_2(G)$, we respectively denote the frattini subgroup and the second center of $G$. The nilpotency class of $G$ is denoted as $cl(G)$, and by $d(G)$ we denote the rank of $G$. By $C_m$ we denote the cyclic group of order $m$ and by $C_m^n$ we denote the $n$ copies of $C_m$. All other unexplained notations are standard. The following results of Curran \cite[Corollaries 3.7, 3.8]{Curran} and Gavioli \cite[Lemma 2]{Gavioli} will be used quite frequently.

\begin{lm}   
Let $G$ be a finite non-abelian $p$-group such that $\mathrm{Aut}_z(G)=Z(\mathrm{Inn}(G))$. Then $Z(G)\le G^{\prime}$ and  $Z(\mathrm{Inn}(G))$ is not cyclic.
\end{lm}

\begin{lm}
Let $G$ be a (not necessarily finite) group with center $Z(G)$ and second center $Z_2(G)$, $A$ an abelian normal subgroup of $G$ with $Z(G)\le A \le Z_2(G)$. Then $A/Z(G)$ embeds in $\mathrm{Hom}(G/A,Z(G))$.
\end{lm}

The next lemma is also of our interest and can be proved using similar arguments as in \cite[Lemma 3]{Alperin}.

\begin{lm}
Let $G$ be any group and $Y$ be a central subgroup of $G$ contained in a normal subgroup $X$ of $G$. Then the group of all automorphisms of $G$ that induce the identity on both $X$ and $G/Y$ is isomorphic to
$\mathrm{Hom}(G/X,Y)$.
\end{lm}

\section{Main Results.}

Let $G/G^{\prime}\simeq\prod_{i=1}^n C_{p^{\alpha_i}}$ and
$Z_2(G)/Z(G)\simeq\prod_{i=1}^m C_{p^{\beta_i}}$ be the decompositions of $G/G'$ and $Z_2(G)/Z(G)$ into cyclic groups,
where  for each $i$, $\alpha_i\ge\alpha_{i+1}$ and  $\beta_i\ge\beta_{i+1}$ are positive integers. Observe that if $\mathrm{Aut}_z(G)=Z(\mathrm{Inn}(G))$, then $Z(G)\le G'$ and hence by Lemma 1.3, $\mathrm{Aut}_z(G)\simeq \mathrm{Hom}(G/G',Z(G))\simeq Z_2(G)/Z(G)$. It thus follows that $d(G)d(Z(G))=d(Z_2(G)/Z(G))$.

\begin{thm}
Let $G$ be a finite $p$-group with cyclic center $Z(G)\simeq C_{p^{\gamma_1}}$. Then $\mathrm{Aut}_z(G)=Z(\mathrm{Inn}(G))$ if and only if   either $G/G'\simeq Z_2(G)/Z(G)$ or $d(G)=d(Z_2(G)/Z(G))$,  $\beta_i =\gamma_1$ for $1\le i\le r$ and $\beta_i =\alpha_i$ for $r+1\le i\le n$, where $r,\;1\le r\le n$, is the largest such that $\alpha_r\ge\gamma_1$.
\end{thm}
{\bf Proof.} First suppose that $\mathrm{Aut}_z(G)=Z(\mathrm{Inn}(G))$. Then  $n=d(G)=d(Z_2(G)/Z(G))=m$ and $\mathrm{exp}(G/G')\ge \mathrm{exp}(Z_2(G)/Z(G))$, because $\mathrm{exp}(Z_2(G)/Z(G))\le \mathrm{exp}(Z(G))$ by \cite[5.2.22]{Rob}. If $\mathrm{exp}(G/G')= \mathrm{exp}(Z_2(G)/Z(G))$, then  
$G/G'\simeq Z_2(G)/Z(G)$; and if $\mathrm{exp}(G/G')> \mathrm{exp}(Z_2(G)/Z(G))$, then 
$\mathrm{exp}(Z_2(G)/Z(G))=\mathrm{exp}(Z(G))$ and hence $\beta_1=\gamma_1$. Suppose $r,\;1\le r\le n$, is the largest such that $\alpha_r\ge\gamma_1$. Then
\[
\mathrm{Aut}_z(G)\simeq \mathrm{Hom}(\textstyle\prod_{i=1}^n C_{p^{\alpha_i}},C_{p^{\gamma_1}})\simeq C_{p^{\gamma_1}}^r\times \textstyle\prod_{i=r+1}^n C_{p^{\alpha_i}}.
\]
But $\mathrm{Aut}_z(G)=  Z_2(G)/Z(G)\simeq \prod_{i=1}^n C_{p^{\beta_i}}$. It follows that $\beta_i =\gamma_1$ for $1\le i\le r$ and $\beta_i =\alpha_i$ for $r+1\le i\le n$.

Conversely, since $Z(G)$ is cyclic, $G$ is purely non-abelian and therefore $|\mathrm{Aut}_z(G)|=|\mathrm{Hom}(G/G^{\prime},Z(G))|$ by \cite[Theorem 1]{Adney}. It is now easy to see, using the hypotheses, that 
 $|\mathrm{Aut}_z(G)|=|Z(\mathrm{Inn}(G))|$ and hence $\mathrm{Aut}_z(G)=Z(\mathrm{Inn}(G))$. 
\hfill$\Box$\\

Let $G$ be a finite $p$-group of order $p^n$ and coclass $\le 4$ such that  $\mathrm{Aut}_z(G)=Z(\mathrm{Inn}(G))$. Then $d(G)d(Z(G))=d(Z_2(G)/Z(G))$ and  $p^{n-5}\le |G'|\le p^{n-2}$. Also, since $Z(\mathrm{Inn}(G))$ cannot be cyclic, $p\le |Z(G)|\le p^3$ and $p^3\le |Z_2(G)|\le p^5$. It is easy to show that $Z(G)$ is cyclic in each case and hence we obtain the following three corollaries.

\begin{cor}
Let $G$ be a finite $p$-group of coclass $2$. Then $\mathrm{Aut}_z(G)=Z(\mathrm{Inn}(G))$ if and only if $Z(G)\simeq C_p$ and $d(G)=d(Z_2(G)/Z(G))=2$.
\end{cor}
{\bf Proof.} First suppose that $\mathrm{Aut}_z(G)=Z(\mathrm{Inn}(G))$. Then $|Z(G)|=p,\;|Z_2(G)|=p^3$ and $Z(\mathrm{Inn}(G))\simeq C_p\times C_p$. It follows by theorem that
either $G/G'\simeq Z_2(G)/Z(G)$ or $d(G)=d(Z_2(G)/Z(G))$. Thus  $d(G)=d(Z_2(G)/Z(G))=2$.
Conversely, suppose that $Z(G)\simeq C_p$ and
$d(G)=d(Z_2(G)/Z(G))=2$. Then $\gamma_1=1$ and $Z_2(G)/Z(G)\simeq C_p\times C_p$. Thus $\alpha_i\ge \gamma_1$ and $\beta_i=\gamma_1=1$ for all $i$. The result now follows by theorem. \hfill $\Box$\\

Using similar arguments, we can prove the following corollary.

\begin{cor}
Let $G$ be a finite $p$-group of coclass $3$. Then $\mathrm{Aut}_z(G)=Z(\mathrm{Inn}(G))$ if and only if either $Z(G)\simeq C_p$ and $d(G)=d(Z_2(G)/Z(G))=2,3$ or $Z(G)\simeq C_{p^2}$ and $Z_2(G)/Z(G)\simeq G/G^{\prime}$.
\end{cor}

\begin{cor}
Let $G$ be a finite $p$-group of coclass $4$. Then $\mathrm{Aut}_z(G)=Z(\mathrm{Inn}(G))$ if and only if one of the followings holds:\\
$(a)$ $Z(G)\simeq C_p$ and $d(G)=d(Z_2(G)/Z(G))=2,3,4$,\\
$(b)$ $Z(G)\simeq C_{p^2}$ and either $(i)$ $Z_2(G)/Z(G)\simeq G/G'$ or $(ii)$ $Z_2(G)/Z(G)\simeq C_{p^2}\times C_p$ and $G/G'\simeq C_{p^3}\times C_p$ or $(iii)$ $Z_2(G)/Z(G)\simeq C_{p^2}\times C_p$ and $G/G'\simeq C_{p^4}\times C_p$,\\
$(c)$ $Z(G)\simeq C_{p^3}$ and $Z_2(G)/Z(G)\simeq G/G^{\prime}$.
\end{cor}
{\bf Proof.} It is not very hard to see that if any of the three conditions hold, then $\mathrm{Aut}_z(G)=Z(\mathrm{Inn}(G))$. Conversely, we prove only part (b), because the other two can be proved using the arguments as in Corollary 2.2. Observe that $p^4\le |Z_2(G)|\le p^5$. If $|Z_2(G)|=p^4$, then $G/G'\simeq Z_2(G)/Z(G)$ by theorem. Next suppose that $|Z_2(G)|=p^5$. If $Z_2(G)/Z(G)\simeq C_p^3$, then  $G/G'\simeq Z_2(G)/Z(G)$ by theorem; and if
$Z_2(G)/Z(G)\simeq C_{p^2}\times C_p$, then $\beta_1=\gamma_1=2$ and $\beta_2<\gamma_1 $, and hence either
 $G/G'\simeq Z_2(G)/Z(G)$ or $G/G'\simeq C_{p^4}\times C_p$ or $G/G'\simeq C_{p^3}\times C_p$. \hfill $\Box$

\section{Groups of order upto $p^7$.} As a consequence of the results proved in last section, we now characterize all finite $p$-groups $G$ of order upto $p^7$ for which $\mathrm{Aut}_z(G)=Z(\mathrm{Inn}(G))$. Proofs are easy, short and even generalize the main results of Sharma and Gumber \cite{Sharma}. The first theorem follows immediately from Corollary 2.2.

\begin{thm}
Let $G$ be a finite $p$-group of order $p^5$ and nilpotency class $3$. Then $\mathrm{Aut}_z(G)=Z(\mathrm{Inn}(G))$ if and only if $Z(G)\simeq C_p$ and $d(G)=d(Z_2(G)/Z(G))=2$.
\end{thm}

\begin{thm}
Let $G$ be a finite $p$-group of order $p^6$ and nilpotency class $3$ or $4$. Then $\mathrm{Aut}_z(G)=Z(\mathrm{Inn}(G))$ if and only if $Z(G)\simeq C_p$ and $d(G)=d(Z_2(G)/Z(G))=2$.
\end{thm}
{\bf Proof.} If $cl(G)=4$, then the result follows by Corollary 2.2. Therefore suppose that $cl(G)=3$. Then  either $Z(G)\simeq C_{p^2}$ and $Z_2(G)/Z(G)\simeq G/G^{\prime}$ or $Z(G)\simeq C_p$ and $d(G)=d(Z_2(G)/Z(G))=2,3$ by Corollary 2.3. We rule out two possibilities to get the result. First suppose that $Z(G)\simeq C_{p^2}$. Then $G/Z(G)$ is a group of order $p^4$ and nilpotency class 2. It follows that $|Z_2(G)|=p^4$ and $|G'|=p^3$, and hence $G/G'$ and $Z_2(G)/Z(G)$ cannot be isomorphic. Next suppose that $|Z(G)|=p$ and $d(G)=d(Z_2(G)/Z(G))=3$. We show that $Z_2(G)$ is abelian. Since $d(G)=cl(G)=3$, $|G'|=p^3$ or $p^2$. If $|G'|=p^3$, then since $G'\le Z(Z_2(G))$, $Z_2(G)$ is abelian. Therefore suppose that $|G'|=p^2$. Since $G/Z_2(G)$ is elementary abelian, $\Phi(G)\le Z_2(G)$ and hence $G'\le C_G(Z_2(G))\le C_G(\Phi(G))$. Thus $\Phi(G)$ is abelian. Let $Z_2(G)=\langle h, \Phi(G)\rangle$. Then, since $h$ centralizes $\Phi(G)$, $Z_2(G)$ is abelian. It now follows by Lemma 1.2 that $Z_2(G)/Z(G)$ embeds in $\mathrm{Hom}(G/Z_2(G), Z(G))\simeq C_p\times C_p$, which is not possible. This completes the proof.  \hfill $\Box$

\begin{thm}
Let $G$ be a finite $p$-group of order $p^7$. Then $\mathrm{Aut}_z(G)=Z(\mathrm{Inn}(G))$ if and only if one of the followings holds:\\
$(i)$ $cl(G)=3, Z(G)\simeq C_p$ and $d(G)=d(Z_2(G)/Z(G))=2,3,4$, \\
$(ii)$ $cl(G)=4$ and either $Z(G)\simeq C_p$,  $d(G)=d(Z_2(G)/Z(G))=2,3$ or $Z(G)\simeq C_{p^2}$ and $Z_2(G)/Z(G)\simeq G/G^{\prime}$,\\
$(iii)$ $cl(G)=5,Z(G)\simeq C_p$ and $d(G)=d(Z_2(G)/Z(G))=2$.
\end{thm}
{\bf Proof.} If $G$ is of nilpotency class 4 (resp. 5), then the result follows from Corollary 2.3 (resp. 2.2). Therefore suppose that $cl(G)=3$. 
Then, by Corollary 2.4, we get the possibilities (a), (b) and (c) for $G$. For the final result, we rule out the  possibilities (b) and (c). First suppose that $Z(G)\simeq C_{p^3}$ and $Z_2(G)/Z(G)\simeq G/G'$. Then, as in above theorem, $G/G'$ is not isomorphic to $Z_2(G)/Z(G)$. Now suppose that $Z(G)\simeq C_{p^2}$ and $Z_2(G)/Z(G)\simeq G/G'\simeq C_p\times C_p$. Then $|G'|=p^5>p^4=|Z_2(G)|$, which is a contradiction to $cl(G)=3$. Next suppose that $Z_2(G)/Z(G)\simeq G/G'\simeq C_p\times C_p\times C_p$. Since $G'\le Z(Z_2(G))$, $Z_2(G)$ is abelian and we get a contradiction by Lemma 1.2. Finally suppose that $Z_2(G)/Z(G)\simeq C_p^2\times C_p$. Then $Z_2(G)=\langle x,y,Z(G)\rangle$, where $x$ and $y$ have orders $p$ and $p^2$ modulo $Z(G)$. Since $x^p\in Z(G)$ and $y\in \Phi(G) =Z_2(G)$, $Z_2(G)$ is abelian and again we get a contradiction by Lemma 1.2. This proves the theorem. \hfill $\Box$

\vspace{.2in} \noindent {\bf Acknowledgement.} The research of the
second author is supported by Council of Scientific and Industrial
Research, Government of India. The same is gratefully acknowledged.


\begin{thebibliography}{20}
\bibitem{Adney} Adney, J. E., Yen, T. (1965).  Automorphisms of a $p$-group. {\em Illinois J. Math.} 9:137-143.
\bibitem{Alperin} Alperin, J. L. (1962). Groups with finitely many automorphisms. {\em Pacific J. Math.}  12:1-5.
\bibitem{Curran} Curran, M. J. (2004).  Finite groups with central automorphism group of minimal order.
{\em Math. Proc. Roy. Irish Acad.}  104 A(2):223-229.
\bibitem{CurMcc} Curran, M. J., McCaughan, D. J. (2001).  Central automorphisms that are almost inner. {\em Comm. Algebra}  29(5):2081-2087. 
\bibitem{Gavioli} Gavioli, N. (1993).  The number of automorphisms of the groups of order $p^7$. {\em Math. Proc. Roy. Irish Acad. Sect. A}  93:177-184. 
\bibitem{Rob} Robinson D. J. S. (1996). {\em A course in the theory of groups}, New York Inc.: Springer-Verlag.
\bibitem{Sharma}  Sharma, M., Gumber, D. (2013).  On central automorphisms of finite
$p$-groups. {\em Comm. Algebra}  41:1117-1122.

 \end{thebibliography}
\end{document}